\documentclass{conm-p-l}

\usepackage{amssymb}

\usepackage{}
\usepackage{hyperref}
\hypersetup{colorlinks=false,pdftitle="Generalized Young Inequality",pdftex}

\newtheorem{thm}{Theorem}[section]
\newtheorem{lem}[thm]{Lemma}
\newtheorem{prop}[thm]{Proposition}
\newtheorem{conj}[thm]{Conjecture}

\theoremstyle{definition}
\newtheorem{definition}[thm]{Definition}

\theoremstyle{remark}
\newtheorem{remark}[thm]{Remark}

\numberwithin{equation}{section}

\newcommand{\bc}{\begin{center}}
\newcommand{\ec}{\end{center}}
\newcommand{\bt}{\begin{tabular}}
\newcommand{\et}{\end{tabular}}
\newcommand{\bea}{\begin{eqnarray}}
\newcommand{\eea}{\end{eqnarray}}
\newcommand{\bean}{\begin{eqnarray*}}
\newcommand{\eean}{\end{eqnarray*}}

\newcommand{\ba}{\begin{array}}
\newcommand{\ea}{\end{array}}

\def\be{\begin{eqnarray}}
\def\ee{\end{eqnarray}}
\def\ben{\begin{eqnarray*}}
\def\een{\end{eqnarray*}}

\newcommand{\ra} {\rightarrow}

\newcommand{\leqa}{\mbox{$ \;\stackrel{(a)}{\leq}\; $}}
\newcommand{\leqb}{\mbox{$ \;\stackrel{(b)}{\leq}\; $}}
\newcommand{\leqc}{\mbox{$ \;\stackrel{(c)}{\leq}\; $}}

\newcommand{\nth}{\frac{1}{n}}

\newcommand{\RL}{{\mathbb R}}

\newcommand{\calC}{\mbox{${\cal C}$}}
\newcommand{\calF}{\mbox{${\cal F}$}}

\newcommand{\Nat}{\mathbb{N}}

\newcommand{\lam}{\lambda}

\def\elabel#1{\label{e:#1}}

\def\sq{$\Box$}

\def\qed{\ifmmode\sq\else{\unskip\nobreak\hfil
\penalty50\hskip1em\null\nobreak\hfil\sq
\parfillskip=0pt\finalhyphendemerits=0\endgraf}\fi\par\medbreak}

\newsavebox{\junk}
\savebox{\junk}[1.6mm]{\hbox{$|\!|\!|$}}

\def\det{{\mathop{\rm det}}}

\def\liminf{\mathop{\rm lim\ inf}}

\def\til={{\widetilde =}}

 \def\eq#1/{(\ref{#1})}

\def\eq#1/{(\ref{e:#1})}

\newcommand{\beqn}[1]{\notes{#1}%
\begin{eqnarray} \elabel{#1}}

\newcommand{\eeqn}{\end{eqnarray} }

\newcommand{\beq}[1]{\notes{#1}%
\begin{equation}\elabel{#1}}

\newcommand{\eeq}{\end{equation}}

\def\bdes{\begin{description}}
\def\edes{\end{description}}

\def\notes#1{}

\newcommand{\calS}{\mathcal{G}}
\newcommand{\calN}{\mathcal{N}}
\newcommand{\calK}{\mathcal{K}}

\renewcommand{\calF}{\mathcal{F}}
\renewcommand{\calC}{\mathcal{C}}

\newcommand{\setS}{\mbox{${\bf s}$}}
\newcommand{\Ys}{Y_{\setS}}
\newcommand{\setT}{\mbox{${\bf t}$}}
\newcommand{\collS}{\calS}

\newcommand{\Mth}{\frac{1}{M}}

\newcommand{\dth}{\frac{1}{d}}
\newcommand{\sumS}{\sum_{\setS\in\collS}}
\newcommand{\prodS}{\prod_{\setS\in\collS}}

\newcommand{\as}{\alpha_{\setS}}

\newcommand{\bs}{\beta_{\setS}}
\newcommand{\gs}{\gamma_{\setS}}

\newcommand{\conv}{\text{conv}}

\begin{document}

\title[Generalizations of Young and Brunn-Minkowski inequalities]{Fractional generalizations of Young and Brunn-Minkowski inequalities}


\author[Bobkov]{Sergey Bobkov}
\address{School of Mathematics\\
University of Minnesota\\
206 Church St. S.E.\\
Minneapolis, MN 55455 USA.
}
\curraddr{}
\email{bobkov@math.umn.edu}
\thanks{}

\author[Madiman]{Mokshay Madiman}
\address{Department of Statistics\\
Yale University\\
24 Hillhouse Avenue\\
New Haven, CT 06511, USA.}
\curraddr{}
\email{mokshay.madiman@yale.edu}
\thanks{}

\author[Wang]{Liyao Wang}
\address{Department of Physics\\
Yale University\\
P.O. Box 208120 \\
New Haven, CT 06520-8120, USA.}
\curraddr{}
\email{liyao.wang@yale.edu}
\thanks{SB was supported in part by NSF grant DMS-0706866, MM was supported
in part by a Junior Faculty Fellowship from Yale University, and LW was supported
in part by the Department of Physics at Yale University.}


\date{}

\begin{abstract}
A generalization of Young's inequality for convolution with sharp constant is conjectured
for scenarios where more than two functions are being convolved,
and it is proven for certain parameter ranges. The conjecture would provide a unified proof
of recent entropy power inequalities of Barron and Madiman, as well as of a (conjectured) generalization of the
Brunn-Minkowski inequality. It is shown that the generalized Brunn-Minkowski conjecture is true for convex sets;
an application of this to the law of large numbers for random sets is described.
\end{abstract}

\maketitle

\section{Introduction}
\label{sec:intro}

Let us denote by $L^p$ the Banach space $L^p(\RL^n,dx)$ of measurable functions
defined on $\RL^n$ whose $p$-th power is integrable with respect to Lebesgue measure $dx$.
In 1912, Young \cite{You12} introduced the fundamental inequality
\be\label{orig-young}
\|f\star g\|_r \leq \|f\|_p \|g\|_q \,,\quad \frac{1}{p}+\frac{1}{q}=\frac{1}{r} +1, \quad 1< p,q,r< +\infty ,
\ee
for functions $f\in L^p$ and $g\in L^q$,
which implies that if two functions are in (possibly different) $L^p$-spaces, then
their convolution is contained in a third $L^p$-space.
In 1972, Leindler \cite{Lei72a} showed the so-called reverse Young inequality,
referring to the fact that the inequality \eqref{orig-young} is reversed when
$0<p,q,r<1$.

For a long time, identification of the best constant that can be put on the right side of
\eqref{orig-young} was an open problem. Eventually, Beckner \cite{Bec75}
proved Young's inequality with the best possible constant. To specify the best
constant, first define $C_p$ by
\be\label{eq:Cp}
C_{p}^{2}=\frac{p^{\frac{1}{p}}}{|p'|^{\frac{1}{p'}}} ,
\ee
where, for any $p\in (0,\infty]$, $p'$ is defined by
\be\label{eq:dual}
\frac{1}{p}+\frac{1}{p'}=1 .
\ee
Note that $p'$ is positive for $p\in (1,\infty)$, and negative for $p\in (0,1)$. Then
the best constant in Young's inequality is $(C_p C_q/ C_{r})^n$.
Soon after, Brascamp and Lieb \cite{BL76b} gave alternative proofs of
both Young's inequality and the reverse Young inequality with this sharp constant;
a simpler and unified proof (of the direct and reverse inequalities) using
transportation arguments was given by Barthe \cite{Bar98a}.
Very recently, an even simpler proof using entropy inequalities was given by Cordero-Erausquin
and Ledoux \cite{CL10}.

Let us remark in passing that a much more general family of inequalities can be proved \cite{BL76b, Lie90};
these are now known as the Brascamp-Lieb inequalities. An optimal transportation proof of these
was given by Barthe \cite{Bar98b}, while an entropy-based proof has recently been
given by Carlen and Cordero-Erausquin \cite{CC09} (cf. Lehec \cite{Leh10:pre}). 
Even more general inequalities are shown using a heat flow
interpolation technique by Bennett, Carbery, Christ and Tao \cite{BCCT10, BCCT08} (see also 
Valdimarsson \cite{Val08, Val10}).
In another direction, Young's inequality can be extended to more general settings than $\RL^n$-- specifically,
to  unimodular locally compact groups (see, e.g., Fournier \cite{Fou77}, Quek and Yap \cite{QY83}, Saeki \cite{Sae90}, 
Baklouti, Smaoui and Ludwig \cite{BSL03} and references therein).

There were several motivations for exploring Young's inequality
with sharp constant-- such as the fact that the optimal constant in
the related Hausdorff-Young inequality (which turns out to be related to the same
$C_p$) gives the definitive formulation of the entropic uncertainty principle,
which is a fundamental result in quantum mechanics. In spite of the
fact that the sharp constant is only very slightly better than 1 for large parameter ranges,
the slight improvement makes all the difference for such applications.

Quite separately from the functional analytic study of $L^p$-norm
inequalities, other mathematical communities were developing inequalities
that would later be seen to be related. Indeed, Brunn, Minkowski  and Lusternik
(cf. \cite{Sch93:book} for the history) developed the famous inequality for volumes of Minkowski
sums that bears their names; this saw enormous development over the following decades,
and became a cornerstone of convex geometry and analysis, apart
from finding numerous applications in a vast variety of fields.
In a completely independent development, Shannon \cite{Sha48}
proposed the so-called ``entropy power inequality'' for entropies
of sums of independent random variables taking values in some Euclidean
space, which was later rigorously proved by Stam \cite{Sta59}.
This inequality in its own way became fundamental in information theory,
emerging as a key tool in proving the so-called converse coding theorems
that show the fundamental limits of various data compression or communication
models.
Subsequently it was noted by several authors that one or both of these inequalities
are related to Young's inequality with sharp constant; indeed proofs of the
Brunn-Minkowski and entropy power inequalities based on Young's inequality
were given by Brascamp-Lieb \cite{BL76a} and Lieb  \cite{Lie78} respectively (see also
Dembo, Cover, Thomas \cite{DCT91}).

Given the history and importance of the results described above, there is clear
intrinsic interest in exploring refinements of them, and in particular
of Young's inequality with sharp constant. While clearly it is impossible
to refine this inequality in the sense of improving the constant, what we explore
in this note is the refinement of it when one is looking at the convolution of {\it more
than two functions}. It turns out that in this case, qualitatively different phenomena
appear that have fascinating connections to random set theory and recent developments in
information theory and probability.

Since we wish to consider $M\geq 2$ functions, let us write $[M]=\{1,2,\ldots,M\}$.
Consider a hypergraph $\collS$ on $[M]$. Recall that a hypergraph
is just a collection of subsets of $[M]$.

Our starting point is the following (unpublished) conjecture made by the
second-named author some years ago; the anonymous referee thought
it might have been discussed before but we have been unable to find a reference.

\begin{conj}\label{conj:young}
Let $\collS$ be a $d$-regular hypergraph on $[M]$.
Let $\{p_{\setS}:\setS\in\collS\}$ and  $r$ be real numbers in $(1,\infty)$ such that
\be\label{eq:p-cond}
\sum_{\setS\in\collS} \frac{1}{p_{\setS}}=|\collS| - \frac{d}{r'} .
\ee
Let $f_{j}, j\in [M]$ be probability density functions on $\RL^n$.
Then
\be\label{eq:main}
\bigg\| \star_{j\in [M]} f_{j} \bigg\|_{r} \leq
\frac{1}{C_{r}^n} \prod_{\setS\in\collS} \bigg[ C_{p_{\setS}}^n
\bigg\| \star_{j\in \setS} f_{j}  \bigg\|_{p_{\setS}} \bigg]^{\frac{1}{d}} .
\ee
Furthermore the inequality is reversed when $\{p_{\setS}:\setS\in\collS\}\cup \{r\} \subset (0,1)$.
\end{conj}

We now outline the main results and organization of this note.
In Section~\ref{sec:young-spl}, Conjecture~\ref{conj:young} is proven for certain parameter ranges.

In Section~\ref{sec:bm}, a conjecture about a generalized Brunn-Minkowski inequality for the Minkowski
sum of more than 2 sets is formulated; it is shown that the conjecture certainly holds for convex sets.
It is also shown in Section~\ref{sec:bm} that Conjecture~\ref{conj:young} implies the conjectured generalized
Brunn-Minkowski inequality for general Borel sets.

In Section~\ref{sec:random}, an application of these generalized Brunn-Minkowski inequalities to the
law of large numbers for random sets is described, after quickly reviewing necessary notions from
the theory of random sets.

Section~\ref{sec:epi} discusses recent generalized entropy power inequalities of 
\cite{MB07}, which gives evidence towards Conjecture~\ref{conj:young} since the former is shown to be a special case of the latter.

Finally, in Section~\ref{sec:sharp-rmk}, we make some remarks on the sharpness of Conjecture~\ref{conj:young}--
in particular, on the question of when extremizers exist.

\section{A special case}
\label{sec:young-spl}

It is appropriate to recall some terminology from discrete
mathematics. A collection $\collS$ of subsets of $[M]$ is called
a {\it hypergraph}, and each set $\setS$ in $\collS$ is called
a {\it hyperedge}. When each hyperedge has cardinality 2,
then $\collS$ can be thought of as the set of edges
of an undirected graph on $m$ labelled vertices.
We interchangeably use ``hypergraph'' and ``collection''
for $\collS$,  ``hyperedge'' and ``set'' for $\setS$ in $\collS$,
and ``vertex'' and ``index'' for $i$ in $[M]$.

The following definitions are  standard.

\begin{definition}
For any index $i$ in $[M]$, define the {\it degree}
of $i$ in $\collS$ as $r(i)=|\{\setT\in\collS: i\in \setT\}|$.

The collection $\collS$ is said to be {\em $d$-regular} if
each index $i$ in $[M]$ has the same degree $d$, i.e.,
if each vertex $i$ appears in exactly $d$ hyperedges of $\collS$.
\end{definition}

The following definition extends the familiar notion of
a partition of a set by allowing
fractional counts. The origin of this notion is unclear
to us, but see 
\cite{SU97:book}.

\begin{definition}
Given a collection $\collS$ of subsets of $[M]$, a function $\gamma:\collS \to [0,1]$,
is called a {\em fractional partition}, if for each $i\in [M]$, we have
$\sum_{\setS\in \collS:i\in \setS} \gs = 1$.
\end{definition}

The following simple lemmas are useful.

\begin{lem}[{\sc Fractional Additivity}]
Let $\{a_{i}:i\in [M]\}$ be an arbitrary collection of real numbers.
For any $\setS\subset [M]$, define $a_{\setS}=\sum_{j\in\setS} a_{j}$.
For any fractional partition $\gamma$ using any hypergraph $\collS$,
$a_{[M]} = \sumS \gs a_{\setS}$.
\end{lem}

\begin{proof}
Interchanging sums implies
\ben\begin{split}
\sumS \gs \sum_{i\in \setS} a_{i}
= \sum_{i\in [M]} a_{i} \sumS \gs {\bf 1}_{\{i\in\setS\}}
= \sum_{i\in [M]} a_{i} .
\end{split}\een
\end{proof}

If the hypergraph $\collS$ is $d$-regular, then
\ben\begin{split}
\sum_{\setS\in\collS, \setS\ni i} \frac{1}{d}
= \sumS \frac{{\bf 1}_{\{i\in\setS\}}}{d}
= 1 ,
\end{split}\een
which motivates the following definition.

\begin{definition}
If $\collS$ is $d$-regular,
$\as=\frac{1}{d}$
defines a fractional partition of $[M]$ using $\collS$,
which we call the {\it degree partition}.
\end{definition}

The following slight extension of H\"older's inequality is useful. We adopt the notation
$f_{\setS}$ for $\prod_{j\in\setS} f_{j}$.

\begin{lem}[Fractional H\"older inequality]
Let $f_{j}, j\in [M]$ be measurable functions on $\RL^{n}$.
Let $\gamma$ be a fractional partition using the hypergraph $\collS$,
and $q_{\setS}$ be  coefficients such that
\be\label{holder-constraint}
\sum_{\setS\in\collS} \frac{\gamma_{\setS}}{q_{\setS}}=\frac{1}{r}.
\ee
Then we have
\be
\| f_{[M]} \|_{r} \leq  \prod_{\setS\in\collS}
\| f_{\setS} \|_{q_{\setS}}^{\gamma_{\setS}} \, .
\ee
\end{lem}

\begin{proof}
Recall that H\"older's inequality says
\ben
\bigg\|\prod_{i\in[M]} f_{i}\bigg\|_{q}\leq \prod_{i\in[M]} \|f_{i}\|_{p_{i}} ,
\een
if $\sum_{i} \frac{1}{p_{i}}=\frac{1}{q}$. (This is traditionally
stated with $q=1$, but it is easy to deduce the form above from that.)
Hence, for any fractional partition $\gamma$ using $\collS$,
\ben
\bigg\| \prod_{j\in[M]} f_{j}\bigg\|_{r}
=\bigg\| \prodS \bigg\{ \prod_{j\in\setS} f_{j} \bigg\}^{\gamma_{\setS}} \bigg\|_{r}
\leq \prodS  \bigg\| \bigg\{ \prod_{j\in\setS} f_{j} \bigg\}^{\gamma_{\setS}} \bigg\|_{p_{\setS}} ,
\een
where $\sumS 1/p_{\setS}=1/r$. But
$\|f_{\setS}^{\gamma_{\setS}}\|_{p_{\setS}} = \|f_{\setS}\|^{\gamma_{\setS}}_{\gamma_{\setS} p_{\setS}}$, so that we obtain the result by setting
$q_{\setS}=\gamma_{\setS} p_{\setS}$ to satisfy the constraint \eqref{holder-constraint}.
\end{proof}

In particular, if $\collS$ be a $d$-regular hypergraph and the
coefficients $q_{\setS}$ satisfy
$\sum_{\setS\in\collS} \frac{1}{q_{\setS}}=\frac{d}{r}$,
we have
\ben
\bigg\| \prod_{j\in [M]} f_{j} \bigg\|_{r} \leq \bigg[ \prod_{\setS\in\collS}
\| f_{\setS} \|_{q_{\setS}} \bigg]^{\frac{1}{d}} .
\een
Combining this elementary observation with the Hausdorff-Young inequality,
Conjecture~\ref{conj:young} follows for a subset of possible parameters.

\begin{thm}\label{thm:subparam}
Conjecture~\ref{conj:young} holds when  $r\geq 2$ and $p_{\setS}\in [1,2]$ for each $\setS$ in $\collS$.
\end{thm}

\begin{proof}
The proof uses the sharp Hausdorff-Young inequality (also
called the Babenko-Beckner inequality). The latter states that if
$f\in L^{p}$ for $p\in [1,2]$, and $\hat{f}$ defined by $\hat{f}(x)= \int e^{2\pi i \langle x,y\rangle} f(y) dy$ is its Fourier
transform, then
\be
\|\hat{f}\|_{p'} \leq C_{p}^n \|f\|_{p} .
\ee
Indeed,
\ben\begin{split}
\big\| \star_{j\in [M]} f_{j} \big\|_{r}
&\leqa C_{r'}^n \big\| \prod_{j\in [M]} \hat{f}_{j} \big\|_{r'}\\
&\leqb C_{r'}^n \prod_{\setS\in\collS}
\bigg[ \big\| \prod_{j\in \setS} \hat{f}_{j}  \big\|_{p'_{\setS}} \bigg]^{\frac{1}{d}} \\
&\leqc C_{r'}^n \prod_{\setS\in\collS}
\bigg[ C_{p_{\setS}}^n \big\| \star_{j\in \setS} f_{j}  \big\|_{p_{\setS}} \bigg]^{\frac{1}{d}} ,
\end{split}\een
where (a) and (c) follow from the Hausdorff-Young inequality, and (b) follows by the fractional H\"older
inequality since \eqref{eq:p-cond} implies that $\sum_{\setS\in\collS} \frac{1}{p'_{\setS}}=\frac{d}{r'}$.
Observing that $C_{r'}=C_{r}^{-1}$ for $r>1$ completes the proof.
\end{proof}

Unfortunately the subset of parameters $p_{\setS}, r$ covered by Theorem~\ref{thm:subparam}
is not the most interesting subset, at least for the applications we have in mind.

\section{Brunn-Minkowski Inequalities}
\label{sec:bm}

Below we always use $|K|$ to denote volume (Lebesgue measure)
of a Borel subset $K$ of Euclidean space of some fixed dimension $n$.
Let $+$ denote the Minkowski sum whenever the addition operation is applied to sets.
Then the classical Brunn-Minkowski inequality states that for
any nonempty Borel sets $K_{1},\ldots, K_{M}$  in $\RL^{n}$,
\ben
|K_{1}+\ldots+K_{M}|^\nth \geq \sum_{j\in [M]} | K_{j} |^\nth .
\een

First we propose the following extended Brunn-Minkowski inequality.

\begin{conj}\label{conj:fracBM}
Let $K_{1},\ldots, K_{M}$ be nonempty Borel sets in $\RL^{n}$.
Then for any fractional partition $\beta$
using the collection $\collS$ of subsets of $[M]$,
\ben
|K_{1}+\ldots+K_{M}|^\nth \geq \sumS \bs \bigg| \sum_{j\in\setS} K_{j} \bigg|^\nth .
\een
\end{conj}

Indeed, observe that this specializes to the usual Brunn-Minkowski inequality when
one takes $\collS$ to be the set of singletons, and each $\bs=1$.

To see the relationship between Young-type and Brunn-Minkowski-type inequalities,
it is useful to define the notion
of R\'enyi entropy, a one-parameter family of entropy-like quantities.
For any random vector $X$ in $\RL^n$ with density $f$, and any $p>1$,
the {\it R\'enyi entropy} of $X$ of order $p$:
\ben
h_p(X)= \frac{p}{p-1}\log \frac{1}{\|f\|_p} ,
\een
where
\ben
\|f\|_p= \bigg(\int_{\RL^n} f^p\,dx\bigg)^{\!1/p}
\een
is the usual $L^p$-norm with respect to Lebesgue measure on $\RL^n$.
The definition of $h_p(X)$ continues to make sense
for $p\in (0,1)$ even though $\|f\|_p$ is then not a norm. There remain the values
$p=0,1,\infty$ on the non-negative half line; for these values, $h_p(X)$ may be defined
``by continuity''. Specifically, as $p\ra 1$, $h_p(X)$ reduces to the Shannon differential entropy
\be\label{def:sha-ent}
h(X)=h_1(X)= -\int_{\RL^n} f(x) \log f(x) dx ,
\ee
and as $p\ra 0$, $h_p(X)$ reduces to
\ben
h_0(X)=  \log |\text{Supp}(f)| ,
\een
where $\text{Supp}(f)$ is the support of the density $f$ (i.e., the closure of the set $\{x\in\RL^n:f(x)>0\}$).

One may also define the  {\it R\'enyi entropy power} of $X$ of order $p$:
\be\label{eq:RenyiEP}
V_p(X)=\exp\bigg\{\frac{2}{n} h_p(X)\bigg\} .
\ee
This reduces to the Shannon entropy power for $p=1$,
and reduces for $p=0$ to
\be\label{eq:0ep}
V_0(X)= |\text{Supp}(f)|^{2/n} .
\ee

\begin{prop}\label{prop:fracBM}
If Conjecture~\ref{conj:young} is true, then Conjecture~\ref{conj:fracBM} is true.
\end{prop}

\begin{proof}
The proof we give is an extension of that used by Dembo, Cover and Thomas \cite{DCT91}
to show that the reverse Young inequality with sharp constant implies the usual Brunn-Minkowski inequality,
and involves taking the limit in an appropriate reformulation of Conjecture~\ref{conj:young} as $r\ra 0$ from above.

Let $X_i$ be random vectors in $\RL^n$ with densities $f_i$ respectively.
The reverse part of Conjecture~\ref{conj:young} asserts that for any $r\in(0,1)$ and $p_{\setS}\in (0,1)$,
\ben
\bigg\| \star_{j\in [M]} f_{j} \bigg\|_{r} \geq
\frac{1}{C_{r}^n} \prod_{\setS\in\collS} \bigg[ C_{p_{\setS}}^n
\bigg\| \star_{j\in \setS} f_{j}  \bigg\|_{p_{\setS}} \bigg]^{\frac{1}{d}} .
\een
Taking the logarithm and rewriting the definition \eqref{eq:RenyiEP} of the R\'enyi entropy power as
$V_p(X)=\|f\|_p^{-2p'/n}$, we have
\be\label{eq:log-conj}\begin{split}
\frac{n}{2r'}\log V_r\bigg(\sum_{i\in[M]} X_i\bigg)
\leq n\log C_{r} &- \frac{n}{d}\sumS \log C_{p_{\setS}}\\
&+ \dth \sumS \frac{n}{2p_{\setS}'} \log V_{p_{\setS}}\bigg(\sum_{i\in\setS} X_i\bigg) .
\end{split}\ee

It is useful to introduce two discrete probability measures $\lam$ and $\kappa$ defined on the hypergraph $\collS$, with
probabilities proportional to $1/p_{\setS}'$ and $1/p_{\setS}$ respectively. 
Let us set 
$L_r=r|\collS|-(r-1)d= r(|\collS|-d/r')$; then the condition \eqref{eq:p-cond},
allows us to write explicitly
\be\label{eq:def-kap}
\kappa_{\setS}  
= \bigg(\frac{r}{L_r}\bigg) \frac{1}{p_{\setS}} \,, \quad \setS\in\collS ,
\ee
and
\be\label{eq:def-lam}
\lambda_{\setS}=\bigg(\frac{r'}{d}\bigg) \frac{1}{p_{\setS}'}   \,, \quad \setS\in\collS ,
\ee
by also using $1/p_{\setS}+1/p_{\setS}' =1$ for the latter.
Then, setting $\Ys=\sum_{i\in\setS} X_i$, \eqref{eq:log-conj} reduces to
\ben
\log V_r(Y_{[M]})
\geq r' \log C_{r}^2 - \frac{r'}{d}\sumS \log C_{p_{\setS}}^2 + \sumS \lambda_{\setS} \log V_{p_{\setS}}(\Ys) .
\een

We wish to write this only in terms of $d, r$, and $\lambda$,
so that we can take $\lambda$ to be fixed and control all other parameters by
tuning $r$ as desired. Towards that end, note that
\ben
r' \log C_{r}^2= -\log |r'|+\frac{r'}{r}\log r
\een
and
\ben\begin{split}
-\frac{r'}{d}\sumS \log C_{p_{\setS}}^2 &= -\frac{r'}{d}\sumS \bigg[ \frac{\log p_{\setS}}{p_{\setS}} - \frac{\log | p_{\setS}' |}{p_{\setS}'} \bigg] \\
&= \sumS \lambda_{\setS} \log | p_{\setS}' | -  \frac{r'}{d} \sumS\log p_{\setS} + \sumS \lambda_{\setS} \log p_{\setS} ,
\end{split}\een
using the definitions \eqref{eq:def-lam} and \eqref{eq:dual} of $\lambda_{\setS}$ and  $p_{\setS}'$. Thus one obtains
\be\label{eq:inter22}\begin{split}
\log V_r(Y_{[M]})
\geq \sumS \lambda_{\setS} &\log V_{p_{\setS}}(\Ys) + \frac{r'}{r}\log r \\
&+ \sumS \lambda_{\setS} \log\bigg[ \frac{ | p_{\setS}' | }{|r'|}\bigg]
 + \sumS \bigg( \lambda_{\setS}-\frac{r'}{d} \bigg)\log p_{\setS} .
\end{split}\ee
The third of the four terms on the right side of \eqref{eq:inter22} simplifies as
\ben
\sumS \lambda_{\setS} \log\bigg[ \frac{ | p_{\setS}' | }{|r'|}\bigg]
=H(\lam)-\log d ,
\een
since 
$|p_{\setS}' | /|r'|=(d\lam_{\setS})^{-1}$ by \eqref{eq:def-lam}, 
where we use $H(\lam)=-\sumS \lambda_{\setS} \log \lambda_{\setS}$ to denote the {\it discrete entropy} of the distribution $\lam$.
Also, the fourth term simplifies as
\ben\begin{split}
\sumS \bigg( \lambda_{\setS}-\frac{r'}{d} \bigg)\log p_{\setS} 
&= \frac{L_{r}}{d(1-r)}  \sumS \kappa_{\setS}  \bigg[\log \frac{1}{\kappa_{\setS}}+\log \frac{r}{L_r} \bigg] \\
&= \frac{L_{r}}{d(1-r)} \bigg[ H(\kappa) +\log \frac{r}{L_r}\bigg],
\end{split}\een
where the first equality follows from the fact that
\ben
\lambda_{\setS}-\frac{r'}{d} = -\frac{r'}{dp_{\setS}}
=-\frac{r'L_{r}\kappa_{\setS}}{dr}
=\frac{L_{r}\kappa_{\setS}}{d(1-r)} 
\een
by successive use of \eqref{eq:def-lam}, \eqref{eq:def-kap} and \eqref{eq:dual},
and from the relation between $p_{\setS}$ and $\kappa_{\setS}$ in \eqref{eq:def-kap}.
With these simplifications \eqref{eq:inter22} can be rewritten as
\be\label{eq:prelimit}\begin{split}
\log V_r(Y_{[M]})
&\geq \sumS \lambda_{\setS} \log V_{p_{\setS}}(\Ys) - \frac{1}{(1-r)}\log r \\
&\quad\quad\quad + H(\lam) -\log d+\frac{L_{r}}{d(1-r)} \bigg[ H(\kappa) +\log \frac{r}{L_r}\bigg] \\
&=\sumS \lambda_{\setS} \log V_{p_{\setS}}(\Ys)  +\frac{L_r-d}{d(1-r)}\log r  \\
&\quad\quad\quad + H(\lam) -\log d+ \frac{L_r}{d(1-r)} \big[ H(\kappa)-\log L_r \big] .
\end{split}\ee

These computations hold for any $r$ and any $\{p_{\setS}\}$, or equivalently,  for any $r$ and any $\lam$.
Let us fix $\lam$; thus one can think of the coefficients $p_{\setS}$
now as functions of $r$. 
We now choose to send $r\downarrow 0$ in \eqref{eq:prelimit}. Then $L_r\ra d$, and
$(L_r-d)\log r= (|\collS|-d)r\log r \ra 0$. Furthermore, from the definitions
\eqref{eq:def-kap} and \eqref{eq:def-lam} of $\kappa$ and $\lambda$,
\ben
\frac{\kappa_{\setS}}{\lam_{\setS}}= \lam_{\setS}^{-1} \frac{r}{L_r} \bigg({1-\frac{1}{p_{\setS}'}} \bigg)
=  \frac{\lam_{\setS}^{-1} r}{L_r} \bigg({1-\frac{d\lambda_{\setS}}{r'}} \bigg)
= \frac{\lam_{\setS}^{-1} r+ (1-r)d}{L_r}
\ra 1 ,
\een
which gives by continuity of the discrete entropy that $H(\kappa)\ra H(\lam)$. Thus, in the limit
as $r\downarrow 0$, the inequality \eqref{eq:prelimit} becomes
\be\label{eq:0limit}
\log V_0(Y_{[M]})\geq \sumS \lambda_{\setS} \log V_{0}(\Ys) + 2[ H(\lam) -\log d ] .
\ee

If  $\text{Supp}(f_i)=K_i$, then $\text{Supp}(\star_{i\in\setS} f_i)=\sum_{i\in\setS} K_i$, which we may denote by $K_{\setS}$; so
\eqref{eq:0limit} simplifies using \eqref{eq:0ep} to
\ben
\nth \log |K_{[M]}| \geq  \sumS \lambda_{\setS} \nth \log |K_{\setS}| +H(\lam) -\log d = \sumS \lambda_{\setS} \log \frac{|K_{\setS}|^\nth}{\lambda_{\setS}} -\log d
\een
The right side is clearly maximized by choosing $\lambda_{\setS}$ proportional to $|K_{\setS}|^\nth$, in which case we obtain
\ben
\log |K_{[M]}|^\nth \geq \log \sumS |K_{\setS}|^\nth -\log d ,
\een
which is precisely the desired result for $d$-regular hypergraphs $\collS$ equipped with the degree partition.
In fact,  assuming the truth of Conjecture~\ref{conj:young}, we have proved that Conjecture~\ref{conj:fracBM} is true
for all regular multihypergraphs (i.e., collections of sets in which a given set may appear
multiple times with different labels, and we keep track of the labels in checking regularity).
The desired result then follows by a bootstrapping argument.
\end{proof}

\begin{remark}
To finish the proof, we used the fact that Conjecture~\ref{conj:fracBM} 
follows from its specialization to $d$-regular multihypergraphs $\collS$ equipped with the degree partition.
While such a bootstrapping capability appears to be folklore in the combinatorics literature, a proof
can be found, e.g., in \cite[Proposition~1]{MG09:isit}. The key point is that the set of all fractional partitions
(when viewed as points in the non-negative orthant of $\RL^{2^{[M]}}$) is a convex, compact set; so linear inequalities hold
for every fractional partition if they hold for every extreme point of the set of fractional partitions. Furthermore, it can be shown
that all these extreme points have rational coordinates, and thus can be viewed as degree partitions
corresponding to certain regular multihypergraphs.
\end{remark}

\begin{remark}
In fact, one can state the following fractional formulation of Conjecture~\ref{conj:young}:
for any fractional partition $\beta$ using the hypergraph $\collS$ on $[M]$,
any  density functions $\{f_{j}, j\in [M]\}$,
and numbers $\{p_{\setS}:\setS\in\collS\}$ and  $r$ satisfying
\be\label{eq:p-cond-frac}
\sum_{\setS\in\collS} \frac{\bs}{p_{\setS}}=\sumS \bs - \frac{1}{r'} ,
\ee
we have
\be\label{eq:main--frac}
\bigg\| \star_{j\in [M]} f_{j} \bigg\|_{r} \leq
\frac{1}{C_{r}^n} \prod_{\setS\in\collS} \bigg[ C_{p_{\setS}}^n
\bigg\| \star_{j\in \setS} f_{j}  \bigg\|_{p_{\setS}} \bigg]^{\bs}
\ee
when $\{p_{\setS}:\setS\in\collS\}\cup \{r\} \subset (1,\infty)$,
and the reverse inequality when $\{p_{\setS}:\setS\in\collS\}\cup \{r\} \subset (0,1)$.
Not surprisingly, this formulation would directly yield Proposition~\ref{prop:fracBM}
via the limiting argument outlined above. (However, although this formulation
appears more general than Conjecture~\ref{conj:young},
they are actually equivalent in keeping with the previous remark.)
\end{remark}

\begin{remark}
Observe that {\it both} the Young and reverse Young inequalities can be compactly expressed in the form \eqref{eq:prelimit},
which holds with the same sign for all positive $p_{\setS}$ and $r$.
\end{remark}

For the special case of convex sets, it is easy to see that Conjecture~\ref{conj:fracBM} is true. The proof relies
on a simple lemma.

\begin{lem}\label{lem:distr-law}
For nonempty convex sets $A$ and $B$, one has the distributive identities
\ben
(a+b)A=aA+bA \quad\text{and}\quad a(A+B)=aA+aB ,
\een
for any non-negative real numbers $a$ and $b$,
whereas these do not hold for general sets.
\end{lem}

\begin{thm}\label{thm:cvx-fracBM}
Let $K_{1},\ldots, K_{M}$ be nonempty convex sets in $\RL^{n}$.
Then for any fractional partition $\beta$
using the collection $\collS$ of subsets of $[M]$,
\be\label{eq:cvx-fracBM}
|K_{1}+\ldots+K_{M}|^\nth \geq \sumS \bs \bigg| \sum_{j\in\setS} K_{j} \bigg|^\nth .
\ee
If the sets $K_j$ are homothetic, one has equality.
\end{thm}

\begin{proof}
By Lemma~\ref{lem:distr-law}, for any fractional partition,
\ben
K_{1}+\ldots+K_{M}=\sumS \bs \sum_{j\in\setS} K_{j} .
\een
Applying the usual Brunn-Minkowski inequality gives
\ben
|K_{1}+\ldots+K_{M}|^\nth \geq \sumS \bigg|\bs \sum_{j\in\setS} K_{j}\bigg|^\nth
= \sumS \bs \bigg|\sum_{j\in\setS} K_{j}\bigg|^\nth .
\een
The equality conditions for the Brunn-Minkowski inequality for convex sets 
require that the sets be homothetic (i.e., equal upto translation and dilatation).
Thus we find that one has equality in \eqref{eq:cvx-fracBM} if and only if the sets
\ben
\bs \sum_{j\in\setS} K_{j} \,,\,\setS\in\collS
\een
are homothetic. This is certainly satisfied if the sets $K_j$ are homothetic.
\end{proof}

Let us note in passing that a different kind of refinement of the Brunn-Minkowski inequality
for convex bodies that captures the ``stability'' of the characterization of extremizers (homothetic
convex bodies) has been recently developed (see, e.g., \cite{FMP09}).

It is interesting to consider adaptations of Theorem~\ref{thm:cvx-fracBM} to 
Gaussian measures. In this context, it is useful to recall the current understanding of
Brunn-Minkowski-type inequalities for Gaussian measure.
The first step towards such an inequality was implicit in Borell's study of log-concave measures \cite{Bor74}; 
in particular, the fact that log-concave measures are characterized by log-concave densities implies that
for Borel sets $K_i\subset \RL^n$, and any $\lam\in [0,1]$,
\be\label{eq:weak-gau}
\gamma(\lam K_{1}+ (1-\lam) K_{2}) \geq \gamma(K_{1})^{\lam} \gamma(K_{2})^{1-\lam} ,
\ee
where $\gamma$ is the standard Gaussian measure on $\RL^n$.
Unlike in the case of Lebesgue measure, however, this log-concavity of measure does
not imply the Gaussian isoperimetric inequality, proved independently by
Sudakov and Tsirelson \cite{ST78} and Borell \cite{Bor75b} (cf. also \cite{Bob97:1}). The latter inequality
asserts that halfspaces are extremal in that they have smallest boundary $\gamma$-measure among
all sets of given $\gamma$-measure. A satisfactory strengthening of \eqref{eq:weak-gau}, which
implies Gaussian isoperimetry, was first obtained by Ehrhard \cite{Ehr83}. In its most general formulation,
due to Borell \cite{Bor03}, it asserts that for Borel sets $K_i\subset \RL^n$ of positive volume, and any $\lam\in [0,1]$,
\be\label{eq:strong-gau}
\Phi^{-1}\circ\gamma(\lam K_{1}+ (1-\lam) K_{2}) \geq \lam \Phi^{-1}\circ\gamma(K_{1}) + (1-\lam) \Phi^{-1}\circ\gamma(K_{2}) ,
\ee
where $\Phi$ is the cumulative distribution function of the one-dimensional standard normal.
(This was proved earlier in \cite{Ehr83} for closed, convex sets, and by Lata{\l}a \cite{Lat96} when one of the sets is
Borel and the other convex.)
The inequality \eqref{eq:strong-gau} has been further generalized by Borell \cite{Bor08}
(cf. Barthe and Huet \cite{BH09} and Gardner and Zvavitch \cite{GZ10}),
where non-convex combinations are also considered.

By an argument very similar to that used in proving Theorem~\ref{thm:cvx-fracBM},
we immediately obtain the following version for Gaussian measure.

\begin{thm}\label{thm:gaus-fracBM}
Let $K_{1},\ldots, K_{M}$ be convex sets of positive volume in $\RL^{n}$.
Suppose  $\beta$ is any  fractional partition
using the collection $\collS$ of subsets of $[M]$, and that the coefficients $\lam_j \geq 0$ satisfy
with $\sum_{j\in[M]} \lam_j =1$. Then we have
\be\label{eq:frac-gau}
\Phi^{-1}\circ\gamma\bigg(\sum_{j\in[M]} \lam_j K_{j}\bigg) \geq \sumS \bs\lam_{\setS} \,\, \Phi^{-1}\circ\gamma\bigg( \sum_{j\in\setS} \frac{\lam_j}{\lam_{\setS}} K_{j} \bigg) ,
\ee
where $\lam_{\setS}=\sum_{i\in\setS} \lam_i$.
\end{thm}

\begin{proof}
Note that
\ben
\sum_{j\in[M]}\lambda_jK_j=\sum_{\setS\in\collS}\bs\lam_{\setS}\sum_{j\in S}\frac{\lambda_j}{\lam_{\setS}}K_j
\een
and
$\sum_{\setS\in\collS}\bs\lam_{\setS}=1$.
Then apply \eqref{eq:strong-gau}.
\end{proof}

Note that the assumption of positive volume (or equivalently positive $\gamma$-measure)
in Theorem~\ref{thm:gaus-fracBM} can be removed, provided we 
adopt the convention $\infty-\infty=-\infty+\infty=-\infty$. For example, if one of the sets is the empty set,
then one should interpret any Minkowski sum of the empty set with any other sets as the empty set,
which would make the right side equal to $-\infty$ and the inequality trivially true. 
It is natural to conjecture that \eqref{eq:frac-gau} continues to hold for all Borel sets.

\section{Applications to random sets}
\label{sec:random}

\subsection{Random sets}

In order to develop the application of Theorem~\ref{thm:cvx-fracBM} to the theory
of random sets, let us first outline some basic features of that theory. We follow the
exposition of Molchanov \cite{Mol93:book}, 
which the reader can consult for more details.

A random closed set is a random element in the space $\calF$ of all closed
subsets (including the empty set $\phi$) of the basic setting space $E=\RL^n$.
To describe the corresponding probability measures, one needs to specify
a topology and $\sigma$-algebra on $\calF$.
For $A\subset \RL^n$, introduce sub-classes of $\calF$ by
\ben
\calF^A=\{F\in\calF:F\cap A=\phi\} \quad,\quad \calF_A=\{F\in\calF:F\cap A\neq\phi\} .
\een
The ``hit-or-miss'' topology  $T_{\calF}$ on the class $\calF$ is the topology generated by
collections of sets of the form
\ben
\calF^K_{G_1,\ldots,G_{M}}=\calF^K \cap \calF_{G_1}\cap \calF_{G_2} \cap \ldots \cap \calF_{G_{M}}  ,
\een
where $K$ runs over the class $\calK$ of compact sets in $\RL^n$, 
and $G_1, \ldots, G_{M}$ lie in the class of  open sets in $\RL^n$.
It is a classical fact that the topological space $(\calF, T_{\calF})$ is compact,
Hausdorff and separable.

A sequence of closed sets $F_{M}, M \geq 1$, converges in $T_{\calF}$ to a certain closed set $F$ if
and only if both the following conditions are valid:
\begin{enumerate}
\item if $K\cap F = \phi$ for a certain compact $K$, then $K\cap F_{M} = \phi$ for all sufficiently large $M$;
\item  if $G \cap F \neq\phi$ for a certain open set $G$, then $G\cap F_{M}\neq \phi$ for all sufficiently large $M$.
\end{enumerate}
We then write $F_{M}\ra_{\calF} F$.

Suppose $\calK$ is the class of compact subsets of $\RL^n$, and 
let $T_{\calK}$ be the topology on $\calK$ induced by $T_{\calF}$. To ensure the convergence of a
sequence $K_{M}, M \geq 1$, of compact sets in $\calK$ an additional condition is required:
there exists a compact $K'$ such that $K_{M} \subset K'$ for all $M\geq 1$.
We then write $K_{M}\ra_\calK K$.

The convergence of compact sets in $\calK$ can be metrized by means of the Hausdorff
metric $\rho_H$ on $\calK$. The Hausdorff distance between two compacts $K_1$ and $K_2$ is defined
as
\ben
\rho_H(K_1, K_2) = \inf\{\epsilon > 0: K_1 \subset K_2^{\epsilon}, \, K_2 \subset K_1^{\epsilon} \} ,
\een
where
$K^{\epsilon} = K+\epsilon B$ is the $\epsilon$-envelope of $K$,
and $B$ denotes the closed ball of unit radius centered at 0.
The Hausdorff distance between two closed sets is defined similarly; however, it can be infinite.

A random closed set is an $\calF$-valued random element, measurable with
respect to the Borel $\sigma$-algebra $\sigma_{\calF}$
generated by $T_{\calF}$ on $\calF$.
Examples of random closed sets include random points and point
processes, random spheres and balls, random half-spaces and hyperplanes etc.
The distribution of a random closed set $A$ is described by the corresponding probability
measure $P$ on $\sigma_{\calF}$, and hence on sets of the type $\calF^K_{G_1,\ldots,G_{M}}$.
Fortunately, $P$ is determined also by its values on $\calF_K$ for $K$ running through $\calK$ only.
In fact, the capacity 
functional of $A$ is defined by
\ben
T_A(K) \,=\, P\{A\in \calF_K \} \,=\, P\{A\cap K\neq \phi\}
\een
for $K\in\calK$.
The properties of $T$ resemble those of the distribution function.

Recall that the support function $s_A$ of a set $A$ is defined by 
\ben
s_{A}(u)=\sup_{x\in A} \,\, \langle u, x\rangle
\een
for any $u\in\RL^n$. Note that if $X$ is a random closed set, 
then its Lebesgue measure or volume $|X|$, its norm
\ben
\|X\| = \sup \{\|x\|: x \in X\},
\een
and its extent in a given direction $s_X(u)$ are usual real-valued random variables.
Also, $\|X\|<\infty$ almost surely if and only if $X$ is compact.

Define $\calC$ to be the class of convex closed sets in $\RL^n$.
A random closed set is said to be convex if its realizations
are almost surely convex, i.e., if $A$ belongs to $\calC$ almost surely.
Similarly, a random compact, convex set is a random closed set
whose realizations lie almost surely in $\calC\cap\calK$.





\subsection{Law of large numbers for random sets}

To formulate a law of large numbers, we first need a notion of expectation for a random set.
Aumann \cite{Aum65} developed such a notion, which was used extensively in the theory of
set-valued functions and related optimization problems; later Artstein and Vitale \cite{AV75}
pioneered its use in the context of random set theory.

We now define the {\it Aumann expectation} of a random compact set $A$.
A random vector $\xi$ in $\RL^n$ (jointly distributed with $A$ on the same
probability space) is said to be a {\it selector} of $A$ if $\xi\in A$ with probability one.
The expectation of $A$ is defined to be the set
\ben
EA = \{E\xi: \xi \,\text{is a selector of }\, A, \, E\xi \,\text{exists} \}.
\een
The condition $E\|A\| < \infty$ is enough to determine that $EA$ is nonempty and compact.
It follows from Aumann \cite{Aum65} that, provided the underlying probability measure is non-atomic,
$EA = E\conv(A)$ and hence $EA$ is convex even for non-convex $A$.
In this case, the expectation $EA$ can also be defined as the convex set having the support
function
\ben
s_{EA}(u)=Es_A(u) \quad,\quad u\in \mathbb{S}^{n-1} ;
\een
this definition continues to make sense for unbounded 
random sets.

The following theorem is due to Vitale \cite{Vit90:1}, and may be considered
a Brunn-Minkowski inequality for random sets.

\begin{thm}\label{thm:Vit}
If $A$ is a random compact set with $E\|A\| < \infty$, then
\ben
|EA|^\nth \geq E|A|^\nth.
\een
\end{thm}


Artstein and Vitale \cite{AV75} developed a law of large numbers for random sets.
Their approach first reduces the general problem to the case of random compact convex sets,
and then proves the result for random compact convex sets by invoking an appropriate
result in the Banach space $C(\mathbb{S}^{n-1})$ and applying it to the support functions of
random sets.

\begin{thm}\label{thm:AV-lln}
Let $A, A_1, A_2,\ldots$ be a sequence of i.i.d. random compact sets with $E\|A\| < \infty$. Then
\ben
\nth \sum_{i=1}^{M} A_i \ra_{\calK} EA  \quad\text{a.s.   as } M\ra\infty .
\een
\end{thm}




The stage is now set for us to state and prove a monotonicity property in the law of large numbers for random sets.

\begin{prop}\label{prop:LLN-mono}
Let $A, A_1, A_2,\ldots$ be a sequence of i.i.d. random compact sets with $E\|A\| < \infty$.
If Conjecture~\ref{conj:fracBM} is true, then
\ben
E\bigg\{\bigg|\Mth \sum_{i=1}^{M} A_i \bigg|^\nth \bigg\}
\een
is a non-decreasing sequence in $M$.
\end{prop}

Since the validity of Conjecture~\ref{conj:fracBM} is known for convex sets,
the statement of Proposition~\ref{prop:LLN-mono}  is also valid for convex sets.
In fact, more is true, but first we need to state a classical result (see, e.g., Beer \cite{Bee74:1}).

\begin{prop}\label{prop:beer}
Suppose $(K_i, i\in\Nat)\subset \calC\cap\calK$, i.e., each $K_i$ is a compact, convex set in $\RL^n$.
If $\rho_H(K_M,K)\rightarrow 0$ as $M\rightarrow\infty$.
and $K\in \calC\cap\calK$, 
then $|K_M|\rightarrow|K|$.
\end{prop}

\begin{thm}\label{thm:LLN-mono}
Let $K, K_1, K_2,\ldots$ be a sequence of i.i.d. random compact {\em convex} sets with $E\|K\| < \infty$.
Then
\ben
E\bigg\{ \bigg|\Mth \sum_{i=1}^{M} K_i \bigg|^\nth \bigg\} \nearrow |EK|^\nth  \quad\text{a.s.   as } M\ra\infty .
\een
In other words, the mean effective radius of the empirical mean based on $M$
observations of the random convex set $A$ is a monotonically non-decreasing sequence (in $M$)
that converges to the effective radius of the Aumann expectation of $K$.
\end{thm}

\begin{proof}
Consider the hypergraph $\collS_{M-1}$ of leave-one-out subsets of $[M]$, i.e.,
$\collS_{M-1}=\{\setS\subset[M]:|\setS|=M-1\}$.
This is a $d$-regular hypergraph with degree $d=M-1$, so Theorem~\ref{thm:cvx-fracBM} implies that
\ben
|K_{1}+\ldots+K_{M}|^\nth \geq \frac{1}{M-1} \sum_{\setS\in\collS_{M-1}} \bigg| \sum_{j\in\setS} K_{j} \bigg|^\nth .
\een
Equivalently,
\be\label{mink-loo}
\bigg|\frac{K_{1}+\ldots+K_{M}}{M}\bigg|^\nth \geq \frac{1}{M} \sum_{\setS\in\collS_{M-1}} \bigg| \frac{\sum_{j\in\setS} K_{j}}{M-1} \bigg|^\nth .
\ee
Setting
\ben
L_{M}=\frac{K_{1}+\ldots+K_{M}}{M} ,
\een
and noting that each of the $M$ summands on the right side of \eqref{mink-loo} has the same law as that of $L_{M-1}$,
we find that
\ben
E\big[ |L_{M}|^\nth \big]
\een
is non-decreasing in $M$.

By Theorem~\ref{thm:Vit}, $E[|L_{M}|^{\frac{1}{n}}]\leq |EL_{M}|^{\frac{1}{n}}$. The i.i.d property and the linearity of the
Aumann expectation yield $EL_{M}=EK$. So $E[|L_{M}|^{\frac{1}{n}}]$ will tend to a
finite limit which is not larger than $|EK|^{\frac{1}{n}}$.
On the other hand, since $L_M\ra_{\calK} EK$ almost surely by Theorem~\ref{thm:AV-lln},
and due to the continuity of the volume functional on $\calK\cap\calC$ asserted by Proposition~\ref{prop:beer},
it follows that the limit of $|L_{M}|$ exists almost surely, and moreover that  
\ben
\lim_{M\ra\infty}  |L_{M}|^{\frac{1}{n}}
=|EK|^{\frac{1}{n}} \quad\text{a.s.}
\een
Then, by Fatou's lemma, one has
\ben
\lim_{M\ra\infty} E[|L_{M}|^{\frac{1}{n}}]
\geq E\big[\liminf_{M\ra\infty} |L_{M}|^{\frac{1}{n}}\big] 
= E\big[\lim_{M\ra\infty}  |L_{M}|^{\frac{1}{n}}\big]
=|EK|^{\frac{1}{n}} \quad\text{a.s.},
\een
which is the desired lower bound.
Combining the bounds yields $\lim_{M\ra\infty} E[|L_{M}|^{\frac{1}{n}}]=|EK|^{\frac{1}{n}}$, and completes the proof.
\end{proof}

\section{Entropy power inequalities}
\label{sec:epi}

We comment here on the connections of Conjecture~\ref{conj:young}
with a recently proved class of so-called entropy power inequalities.

For a $\RL^{n}$-valued random vector $X$ with density $f$
with respect to the Lebesgue measure  on $\RL^n$,
the entropy (sometimes called differential entropy or Boltzmann--Shannon entropy) is given by \eqref{def:sha-ent},
and the (Shannon) entropy power of $X$ is
$\calN(X)=e^{2h(X)/n}$.
We limit ourselves to random vectors $X$
with $h(X)<+\infty$; in this case, $\calN(X)$ 
is a non-negative real number.

Building on work of 
\cite{MB07} and resolving a conjecture they made, 
\cite{MG09:isit} recently showed the following result.

\begin{thm}\label{thm:epi}
Let $X_{1},\ldots, X_{M}$ be independent $\RL^{n}$-valued random vectors,
such that the entropy of each exists and is finite.
Let $\beta$ be a fractional partition using a collection $\collS$ of subsets of $[M]$. Then
\ben
\calN(X_{1}+\ldots+X_{M}) \geq \sum_{\setS\in\collS}
\bs \calN\bigg(\sum_{j\in\setS} X_{j}\bigg) .
\een
Equality holds if 
all the $X_{i}$ are normal with proportional
covariance matrices.
\end{thm}

Let us briefly mention some specializations of Theorem~\ref{thm:epi}.
If $\calS$ is an arbitrary hypergraph on $[M]$, \cite{MB07} showed that
\be\label{mb-epi}
\calN(X_{1}+\ldots+X_{M}) \geq \frac{1}{d} \sum_{\setS\in\calS} \calN \big( \sum_{j\in\setS} X_{j}\big)  ,
\ee
where $d$ is the maximum number of hyperedges in $\calS$ in which any one vertex appears (and in particular
for $d$-regular hypergraphs).
Choosing $\calS$ to be the class $\collS_{M-1}$ of all sets of $M-1$ elements
yields $d=M-1$ and hence
\be\label{abbn-epi}
\calN(X_{1}+\ldots+X_{M}) \geq \frac{1}{M-1} \sum_{i\in [M]} \calN\big( \sum_{j\neq i} X_{j}\big) .
\ee
This inequality was proved by Artstein, Ball, Barthe and Naor \cite{ABBN04:1}, and was used by them
to affirmatively resolve the long-standing conjecture of monotonicity in Barron's
entropic central limit theorem \cite{Bar86}.
Choosing $\calS$ to be the class $\collS_{1}$ of all singletons in \eqref{mb-epi}
yields $d=1$ and hence
\be\label{sha-epi}
\calN(X_{1}+\ldots+X_{M}) \geq \sum_{j\in [M]} \calN(X_{j}) ,
\ee
which is the classical Shannon-Stam entropy power inequality \cite{Sha48, Sta59}.
This is already a nontrivial and interesting inequality, implying (as implicitly contained in \cite{Sta59}) 
for instance the logarithmic Sobolev inequality for the Gaussian usually attributed to Gross \cite{Gro75}.

Theorem~\ref{thm:epi} is related to Conjecture~\ref{conj:young}; indeed the former
follows from the latter and thus provides some evidence towards the validity of
Conjecture~\ref{conj:young}. The proof of this implication is very similar to that
of Proposition~\ref{prop:fracBM}, except that one takes the limit $r\ra 1$ instead of $r\ra 0$
in the form \eqref{eq:prelimit} of Conjecture~\ref{conj:young}.

\section{Remarks on the sharpness of Conjecture~\ref{conj:young}}
\label{sec:sharp-rmk}

Consider the following simple case of Conjecture~\ref{conj:young} (we only
consider the generalization of Young's inequality, although similar comments
can be made about reverse Young), corresponding to $n=1$, $M=3$ and $d=2$.
If  $\frac{1}{p }+\frac{1}{q}+\frac{1}{t}=3-\frac{2}{r'}$, then
\ben
{\|}f_1{\star}f_2{\star}f_3{\|}_r
\leq  C_{r'}\sqrt{C_pC_qC_t} \,
{{\|}f_1{\star}f_2{\|}}^{\frac{1}{2}}_p \, {{\|}f_2{\star}f_3{\|}}^{\frac{1}{2}}_q \, {{\|}f_3{\star}f_1{\|}}^{\frac{1}{2}}_t .
\een
Given that Young's inequality with sharp constant is (of course!) sharp, and that equality can only be attained
for Gaussians, it is natural to expect that a similar fact holds for Conjecture~\ref{conj:young}. However,
it turns out that this is not quite the case.

Take $f_i$ to be the density of the non-degenerate normal distribution $N(\mu_i,\sigma_i^2)$ 
with mean $\mu_i$ and variance $\sigma_i^2$, and plug them into the above inequality to get
\be\label{gaus-ineq}
x^{\frac{1}{4p}}y^{\frac{1}{4q}}(2-x-y)^{\frac{1}{4t}}\leqslant {\bigg(\frac{r'}{p'}\bigg)}^{\frac{1}{4p}}
{\bigg(\frac{r'}{q'}\bigg)}^{\frac{1}{4q}} {\bigg(\frac{r'}{t'}\bigg)}^{\frac{1}{4t}}
\ee
where
\ben
x=\frac{\sigma_2^2+{\sigma_3}^2}{\sigma_1^2+{\sigma_2}^2+{\sigma_3}^2} \quad\text{and}\quad
y=\frac{\sigma_1^2+{\sigma_3}^2}{\sigma_1^2+{\sigma_2}^2+{\sigma_3}^2} \, .
\een
Note that $(x,y)$ lies in the region $\{(x,y)\in \textbf{R}^2|x<1,y<1,x+y>1\}$.
Simple calculus shows that if
\be\label{cond-tight}
r'< \min\{p', q', t'\} ,
\ee
then \eqref{gaus-ineq} is sharp.
On the other hand, if the condition \eqref{cond-tight} is violated,
then the right side of \eqref{gaus-ineq} still bounds the left side from above, but
it is not the best bound for the function on the left side. In the rest of this section,
we make some remarks that attempt to shed light on this observation, which is
somewhat unexpected in view of the fact that consequences of Conjecture~\ref{conj:young}
such as Conjecture~\ref{conj:fracBM} and Theorem~\ref{thm:epi} are clearly tight (for homothetic
convex bodies and Gaussians with proportional covariance matrices respectively).

Let us first examine the way in which Conjecture~\ref{conj:young}
implies Conjecture~\ref{conj:fracBM} and Theorem~\ref{thm:epi}.
The strategy was to let $r$ go to some limit (either 0 or 1), while keeping
the coefficients $\lam$ constant. This yielded a limit inequality for any
fixed $\lam$, which was then optimized over $\lam$ to obtain the desired
conclusion. Furthermore, in both the Brunn-Minkowski and entropy power
contexts, the optimal choice of $\lam$ happens to be such that each $\lam_{\setS}$ is always
bounded from above by $1/d$ (or in other words, $r'<\min\{p_{\setS}':\setS\in\collS\}$,
which is condition \eqref{cond-tight} for the general case).
Thus the source of the looseness appears to lie
in the fact that there is an optimization of the inequality that
has not been performed.

Note that the optimal choice of $\lam$ in the preceding discussion {\it depends on the functions $f_i$}.
This suggests that it may be interesting to consider the following problem: Fix all the functions $f_i$,
as well as the parameter $r$ and the $d$-regular hypergraph $\collS$,
 in Conjecture~\ref{conj:young}. Assuming that the conjecture is true,
 what are the best constants $\{p_{\setS}|\setS\in\collS\}$ such that the inequality will hold? (In other
 words, what is the optimized form of the conjectured inequality
{\it without} taking a limit in $r$?) Furthermore, does such an optimization always yield a tight bound on
 the left side of the conjectured inequality, which is achieved for Gaussians?

While we are not able to completely answer these questions, we give some indications.
Using the reformulation \eqref{eq:prelimit} of Conjecture~\ref{conj:young} in terms of R\'{e}nyi entropy powers,
our goal is now  to maximize the right side of \eqref{eq:prelimit} over choice of $\lambda$
(which determines $\{p_{\setS}|\setS\in\collS\}$ and $\kappa$), for fixed functions.
The following simple lemma is useful.

\begin{lem}\label{lem:opt}
Define $\varphi_f(p)=\log\int_X|f|^pd\mu$,where $\mu$ is any measure on the measure space X.
Let $E=\{0<p<\infty|\varphi_f(p)<\infty\}$.
Then
\begin{enumerate}
\item $E$ is a convex set.
\item On $E$, $\varphi_f(p)$ is a convex function in $p$.
\item On $E$, $\varphi_f(p)$ is a continuous function.
\item In the interior of $E$, $\varphi_f(p)$ is infinitely differentiable.
\item In the interior of $E$, $\varphi_f(p)-p\frac{d\varphi_f(p)}{dp}$ is a non-increasing function of $p$.
\end{enumerate}
\end{lem}

\begin{proof}
The first two parts are classical-- indeed, the second is Lyapunov's inequality.
Continuity of $\varphi_f(p)$  on the interior of $E$ follows from its convexity.
Also, $E$ must be an interval since it is convex-- if it includes an endpoint,
use dominated convergence to show that it is left (respectively, right) continuous at the
right (respectively, left) endpoint.

For part (4), suppose $(p_1,p_2)$ is a subset of the interior of $E$, so that
$p_2+\epsilon\in E$ and $p_1-\epsilon\in E$ for some $\epsilon>0$.
Let $p,q\in (p_1,p_2)$, with $p$ fixed, and $q\neq p$.
Note that on $\{|f|>0\}$, $\frac{|f|^{q}-|f|^p}{q-p}=\log(|f|)|f|^\xi$,
where $\xi$ is between $p$ and $q$. 
On the set $\{0<|f|\leqslant1\}$, bound $|\log|f||$ by $M_1|f|^{-\epsilon}$,
$|f|^\xi$ by $|f|^{p_1}$ and $|\frac{|f|^{q}-|f|^p}{q-p}|$ by $M_1|f|^{p_1-\epsilon}$. 
On the set $\{|f|\geqslant1\}$,
bound $|\log|f||$ by $M_2|f|^{\epsilon}$, $|f|^\xi$ by $|f|^{p_2}$ and $|\frac{|f|^{q}-|f|^p}{q-p}|$
by $M_2|f|^{p_2+\epsilon}$. Now use dominated convergence to get the desired result. Similarly for higher order derivatives.

For part (5), note that in the interior of $E$, the derivative of $\varphi_f(p)-p\frac{d\varphi_f(p)}{dp}$ is simply $-p\varphi''_f(p)$,
which is smaller than or equal to zero due to part (2). So it is a non-increasing function in the interior of $E$. 
\end{proof}

We can now apply the Lagrange multiplier method to obtain a necessary condition for 
optimal coefficients $\{p_{\setS}|\setS\in\collS\}$. (It is not known to be sufficient since
the objective function does not appear to be concave.) In the following proposition, we
adopt the notation $f^*_{\setS}=\star_{i\in\setS} f_i$, and use the fact that 
for any density $f$, the quantity $\varphi_f(p)-p\frac{d\varphi_f(p)}{dp}$
from Lemma~\ref{lem:opt} can also be written in terms of the entropy of the new density
function $\frac{{f}^{p}}{\int {f}^{p}dx}$.

\begin{prop}
Suppose $f_i$ are densities on $\RL^n$ such that
$\|f^*_{\setS}\|_{p}$ is finite for all $p\in (0,+\infty)$ and each $\setS\in\collS$.
Then if the set of nonnegative real values
$\{p_{\setS}|\setS\in\collS\}$ maximizes the right side of \eqref{eq:prelimit},
there must exist a constant $\beta\in\RL$ such that the stationary conditions
\ben
\log\bigg[\frac{|1-p_{\setS}|}{{p_{\setS}}^2}\bigg]=\frac{2h(F_{\setS})}{n}+\beta \quad\text{for all $\setS\in\collS$}
\een
and
\ben
\sumS \frac{1}{p_{\setS}} =|\mathcal{\collS}|-\frac{d}{r'}
\een
hold, where
\ben
F_{\setS}= \frac{{(f^*_{\setS})}^{p_{\setS}}}{\int_{\mathbb{R}^n}{(f^*_{\setS})}^{p_{\setS}}dx} .
\een
\end{prop}



The equations above do not seem to be explicitly solvable in general.
However, when each $f_j$ is a centered non-degenerate Gaussian with covariance matrix $K_j$,
the system of equations above becomes explicitly solvable.
Moreover, if one substitutes these values
of $p_{\setS}$ into the right side of \eqref{eq:prelimit}, one obtain the inequality
\ben
\det^\nth\bigg(\sum_{j\in [M]}K_j\bigg) \geq \dth \sumS \det^{\nth} \bigg(\sum_{j\in \setS}K_j\bigg) 
\een
by tedious but entirely elementary calculations. 
Observe that this is a special case of both Theorem~\ref{thm:epi} (applied to Gaussians)
and Theorem~\ref{thm:cvx-fracBM} (applied to ellipsoids), and that it is tight-- in particular, 
it holds with equality if the covariance matrices $K_j$ are proportional.

\vspace{.05in}
\noindent{\bf Acknowledgments.}
We are grateful to Professor Richard Vitale for help with references,
and an anonymous referee for pointing out a number of typos and helping to
improve exposition.

\bibliographystyle{amsplain}
%
%
\providecommand{\bysame}{\leavevmode\hbox to3em{\hrulefill}\thinspace}
\providecommand{\MR}{\relax\ifhmode\unskip\space\fi MR }
\providecommand{\MRhref}[2]{%
  \href{http://www.ams.org/mathscinet-getitem?mr=#1}{#2}
}
\providecommand{\href}[2]{#2}

\end{document}